\RequirePackage[british]{babel}
\documentclass[reqno,a4paper,12pt]{amsart}
\usepackage{a4wide}
\usepackage[utf8x]{inputenc}
\usepackage{amsmath,amssymb,amstext,amsthm,amscd,mathrsfs,eucal}
\usepackage{graphicx,color}
\usepackage{array}
\usepackage{hyperref}
\hypersetup{%
  pdftitle   = {Metric Lie n-algebras and double extensions},
  pdfkeywords = {Lie, Filippov, n-algebra, metric, double extension},
  pdfauthor  = {José Figueroa-O'Farrill},
  pdfcreator = {\LaTeX\ with package \flqq hyperref\frqq}
}
\PrerenderUnicode{éÉ}
\newcommand{\fgl}{\mathfrak{gl}}

\newcommand{\RR}{\mathbb{R}}

\newcommand{\be}{\boldsymbol{e}}

\newcommand{\Cbar}{\overline{C}}
\newcommand{\Ubar}{\overline{U}}
\newcommand{\Wbar}{\overline{W}}

\DeclareMathOperator{\End}{End}
\DeclareMathOperator{\Der}{Der}
\DeclareMathOperator{\Rad}{Rad}

\DeclareMathOperator{\ad}{ad}

\theoremstyle{plain}
\newtheorem{lemma}{Lemma}
\newtheorem{proposition}[lemma]{Proposition}
\newtheorem{theorem}[lemma]{Theorem}
\newtheorem{corollary}[lemma]{Corollary}
\theoremstyle{definition}
\newtheorem{definition}[lemma]{Definition}

\newcommand{\MUNCH}[1]{\relax}

\allowdisplaybreaks
\begin{document}
\title{Metric Lie n-algebras and double extensions}
\author{José Figueroa-O'Farrill}
\address{Maxwell Institute and School of Mathematics, University of
  Edinburgh, UK}
\email{J.M.Figueroa@ed.ac.uk}
\date{\today}
\begin{abstract}
  We prove a structure theorem for Lie $n$-algebras possessing an
  invariant inner product.  We define the notion of a double extension
  of a metric Lie $n$-algebra by another Lie $n$-algebra and prove
  that all metric Lie $n$-algebras are obtained from the simple and
  one-dimensional ones by iterating the operations of orthogonal
  direct sum and double extension.
\end{abstract}
\maketitle
\tableofcontents

\section{Introduction}
\label{sec:intro}

A (finite-dimensional, real) Lie $n$-algebra consists of a
finite-dimensional real vector space $V$ together with a linear map
$\Phi: \Lambda^n V \to V$, denoted simply as an $n$-bracket, obeying a
generalisation of the Jacobi identity.  To define it, let us recall
that an endomorphism $D\in \End V$ is said to be a \textbf{derivation}
if
\begin{equation*}
  D[x_1\dots x_n] = [Dx_1\dots x_n] + \dots + [x_1 \dots Dx_n]~,
\end{equation*}
for all $x_i \in V$.  Then $(V,\Phi)$ defines a \textbf{Lie
  $n$-algebra} if the endomorphisms $\ad_{x_1\dots x_{n-1}} \in \End
V$, defined by $\ad_{x_1 \dots x_{n-1}} y = [x_1 \dots x_{n-1} y]$,
are derivations.  When $n=2$ this clearly agrees with the Jacobi
identity of a Lie algebra.  For $n>2$ we will call it the $n$-Jacobi
identity.  The vector space of derivations is a Lie subalgebra of
$\fgl(V)$ denoted $\Der V$.  The derivations $\ad_{x_1 \dots x_{n-1}}
\in \Der V$ span the ideal $\ad V \lhd \Der V$ consisting of
\textbf{inner derivations}.

From now on, whenever we write Lie $n$-algebra, we will assume that
$n>2$ unless otherwise stated.  In this paper we will only work with
finite-dimensional real Lie $n$-algebras.

Lie $n$-algebras were introduced by Filippov \cite{Filippov} and have
been studied further by a number of people.  We mention here only two
outstanding works beyond Filippov's original paper: the pioneering
work of Kasymov \cite{Kasymov} and the PhD thesis of Ling
\cite{LingSimple}.  Kasymov studied the various notions of solvability
and nilpotency for Lie $n$-algebras, introduced the notion of
representation of a Lie $n$-algebra and proved an Engel-type theorem
and a Cartan-like criterion for solvability.  Ling classified simple
Lie $n$-algebras and proved a very useful Levi-type decomposition.  It
is perhaps remarkable that most structural results in the theory of
Lie $n$-algebras are actually consequences of similar results for the
Lie algebra of derivations.  In this sense it is to be expected that
results for Lie algebras should have their analogue in the theory of
Lie $n$-algebras; although it seems that Lie $n$-algebras become more
and more rare as $n$ increases, due perhaps to the fact that as $n$
increases, the $n$-Jacobi identity imposes more and more conditions.

For example, over the complex numbers there is up to isomorphism a
unique simple Lie $n$-algebra for every $n>2$, of dimension $n+1$ and
whose $n$-bracket is given relative to a basis $(\be_i)$ by
\begin{equation*}
  [\be_1\dots\widehat{\be_i}\dots\be_{n+1}] = (-1)^i \be_i~,
\end{equation*}
where a hat denotes omission.  Over the reals, they are all given by
attaching a sign $\varepsilon_i$ to each $\be_i$ on the right-hand
side of the bracket.

A class of Lie $n$-algebras which have appeared naturally in
mathematical physics are those which possess a nondegenerate inner
product which is invariant under the inner derivations.  We call them
\textbf{metric Lie $n$-algebras}.  They seem to have arisen for the
first time in work of Papadopoulos and the author \cite{FOPPluecker}
in the classification of maximally supersymmetric type IIB
supergravity backgrounds \cite{FOPMax}, and more recently, for the
case of $n=3$, in the work of Bagger and Lambert \cite{BL1,BL2} and
Gustavsson \cite{GustavssonAlgM2} on a superconformal field theory for
multiple M2-branes.  It is this latter work which has revived the
interest of part of the mathematical physics community on metric Lie
$n$-algebras.

Metric Lie algebras are not as well understood as the simple Lie
algebras; although, shy of a classification, a number of structural
results are known.  It is a classic result that Lie algebras
possessing a positive-definite invariant inner product are reductive,
whence isomorphic to an orthogonal direct sum of simple and
one-dimensional Lie algebras.  In lorentzian signature (i.e., index 1)
there is a classification due to Medina \cite{MedinaLorentzian}.  The
indecomposable lorentzian Lie algebras are constructed out of the
one-dimensional Lie algebra by iterating two constructions: orthogonal
direct sum and \emph{double extension}.  This was later extended by
Medina and Revoy \cite{MedinaRevoy} (see also work of Stanciu and the
author \cite{FSalgebra}), who showed that indecomposable metric Lie
algebras are constructed by again iterating the operations of direct
sum and the (generalised) double extension, using again as ingredients
the simple and one-dimensional Lie algebras.  This was used in
\cite{MedinaLorentzian} to construct all possible indecomposable
metric Lie algebras of index 2 (i.e., signature $(2,p)$).  Contrary to
the lorentzian case, there is a certain ambiguity in this
construction, which prompted Kath and Olbrich \cite{KathOlbrich2p} to
approach the classification problem for metric Lie algebras from a
cohomological perspective.  In particular they classified
indecomposable metric Lie algebras with index 2, a result which had
been announced in \cite{baumkath2p}.  For more indefinite signatures,
the classification problem is still largely open.

Much less is known about metric Lie $n$-algebras.  There is a
classification for euclidean \cite{NagykLie} (see also \cite{GPkLie})
and lorentzian \cite{JMFLorNLie} metric Lie $n$-algebras and also a
classification of index-2 metric Lie 3-algebras \cite{2p3Lie}.  In
that paper there is also a structure theorem for metric Lie 3-algebras
and a definition of double extension.  In this note we will extend
these results to $n>3$.  We prove a structure theorem for metric Lie
$n$-algebras and in particular introduce the notion of a double
extension of a metric Lie $n$-algebra by another Lie $n$-algebra.

\section*{Acknowledgments}

It is a pleasure to thank Paul de Medeiros and Elena Méndez-Escobar
for many entertaining and illuminating $n$-algebraic discussions.

\section{Metric Lie $n$-algebras}
\label{sec:metric-n-lie}

We recall that a metric Lie $n$-algebra is a triple $(V,\Phi,b)$
consisting of a finite-dimensional real vector space $V$, a linear map
$\Phi: \Lambda^n V \to V$, denoted simply by an $n$-bracket, and a
nondegenerate symmetric bilinear form $b:S^2 V \to \RR$, denoted
simply by $\left<-,-\right>$, subject to the $n$-Jacobi identity
\begin{equation}
  \label{eq:n-Jacobi}
  [x_1 \dots x_{n-1} [y_1 \dots y_n]] =
  [[x_1 \dots x_{n-1} y_1] \dots y_n] + \dots +
  [y_1 \dots [x_1 \dots x_{n-1} y_n]]~,
\end{equation}
and the invariance condition of the inner product
\begin{equation}
  \label{eq:adinvariance}
  \left<[x_1 \dots x_{n-1} y_1], y_2\right> = -
  \left<[x_1 \dots x_{n-1} y_2], y_1\right>~,
\end{equation}
for all $x_i, y_i \in V$.

Given two metric Lie $n$-algebras $(V_1,\Phi_1,b_1)$ and
$(V_2,\Phi_2,b_2)$, we may form their \textbf{orthogonal direct sum}
$(V_1\oplus V_2,\Phi_1\oplus \Phi_2, b_1 \oplus b_2)$, by declaring
that
\begin{align*}
  [x_1 x_2 y_1 \dots y_{n-2}] = 0 \qquad\text{and}\qquad
  \left<x_1,x_2\right> = 0~,
\end{align*}
for all $x_i\in V_i$ and all $y_i\in V_1 \oplus V_2$.  The resulting
object is again a metric Lie $n$-algebra.  A metric Lie $n$-algebra is
said to be \textbf{indecomposable} if it is not isomorphic to an
orthogonal direct sum of metric Lie $n$-algebras $(V_1\oplus V_2,
\Phi_1\oplus \Phi_2, b_1\oplus b_2)$ with $\dim V_i > 0$.  In order to
classify the metric Lie $n$-algebras, it is clearly enough to classify
the indecomposable ones.  In Section \ref{sec:structure} we will prove
a structure theorem for indecomposable Lie $n$-algebras.

\subsection{Basic facts about Lie $n$-algebras}
\label{sec:basic}

From now on let $(V,\Phi)$ be a Lie $n$-algebra.  Given subspaces $W_i
\subset V$, we will let
\begin{equation*}
  [W_1\dots W_n] = \left\{ [w_1\dots w_n] \middle | w_i \in
    W_i\right\}~.
\end{equation*}

We will use freely the notions of subalgebra, ideal and homomorphisms
as reviewed in \cite{JMFLorNLie}.  In particular a \textbf{subalgebra}
$W < V$ is a subspace $W \subset V$ such that $[W\dots W] \subset W$,
whereas an \textbf{ideal} $I \lhd V$ is a subspace $I \subset V$ such
that $[I V \dots V] \subset I$.  A linear map $\phi: V_1 \to V_2$
between Lie $n$-algebras is a \textbf{homomorphism} if $\phi[x_1 \dots
x_n] = [\phi(x_1) \dots\phi(x_n)]$, for all $x_i\in V_1$.  An
\textbf{isomorphism} is a bijective homomorphism.  There is a
one-to-one correspondence between ideals and homomorphisms and all the
standard theorems hold.  In particular, intersection and sums of
ideals are ideals.  An ideal $I\lhd V$ is said to be \textbf{minimal}
if any other ideal $J\lhd V$ contained in $I$ is either $0$ or $I$.
Dually, an ideal $I\lhd V$ is said to be \textbf{maximal} if any other
ideal $J\lhd V$ containing $I$ is either $V$ and $I$.  If $I \lhd V$
is any ideal, we define the \textbf{centraliser} $Z(I)$ of $I$ to be
the subalgebra defined by $[Z(I) I V \dots V] = 0$. Taking $V$ as an
ideal of itself, we define the \textbf{centre} $Z(V)$ by the condition
$[Z(V) V \dots V]=0$.  A Lie $n$-algebra is said to be \textbf{simple}
if it has no proper ideals and $\dim V > 1$.

\begin{lemma}\label{le:simplequot}
  If $I\lhd V$ is a maximal ideal, then $V/I$ is simple or
  one-dimensional.
\end{lemma}

Simple Lie $n$-algebras have been classified.

\begin{theorem}[{\cite[§3]{LingSimple}}]\label{th:simple}
  A simple real Lie $n$-algebra is isomorphic to one of the
  ($n+1$)-dimensional Lie $n$-algebras defined, relative to a basis
  $\be_i$, by
  \begin{equation}
    \label{eq:simple-n-Lie}
    [\be_1 \dots  \widehat{\be_i} \dots  \be_{n+1}] = (-1)^i \varepsilon_i
    \be_i~,
  \end{equation}
  where a hat denotes omission and where the $\varepsilon_i$ are
  signs.
\end{theorem}

It is plain to see that simple real Lie $n$-algebras admit invariant
inner products of any signature.  Indeed, the Lie $n$-algebra in
\eqref{eq:simple-n-Lie} leaves invariant the diagonal inner product
with entries $(\varepsilon_1, \dots, \varepsilon_{n+1})$.

Complementary to the notion of semisimplicity is that of solvability.
As shown by Kasymov \cite{Kasymov}, there is a whole spectrum of
notions of solvability for Lie $n$-algebras.  However we will use here
the original notion introduced by Filippov \cite{Filippov}.  Let $I
\lhd V$ be an ideal.  We define inductively a sequence of ideals
\begin{equation}
  \label{eq:sequence}
  I^{(0)} = I \qquad\text{and}\qquad
  I^{(k+1)} = [I^{(k)} \dots I^{(k)}] \subset I^{(k)}~.
\end{equation}
We say that $I$ is \textbf{solvable} if $I^{(s)}=0$ for some $s$, and
we say that $V$ is \textbf{solvable} if it is solvable as an ideal of
itself.  If $I,J\lhd V$ are solvable ideals, so is their sum $I + J$,
leading to the notion of a maximal solvable ideal $\Rad V$, known as
the \textbf{radical} of $V$.  A Lie $n$-algebra $V$ is said to be
\textbf{semisimple} if $\Rad V=0$.  Ling \cite{LingSimple} showed that
a semisimple Lie $n$-algebra is isomorphic to the direct sum of its
simple ideals.  The following result is due to Filippov \cite{Filippov}
and can be paraphrased as saying that the radical is a
\emph{characteristic} ideal.

\begin{theorem}[{\cite[Theorem~1]{Filippov}}]\label{th:RadChar}
  Let $V$ be a Lie $n$-algebra.  Then $D\Rad V \subset \Rad V$ for
  every derivation $D \in \Der V$.
\end{theorem}

We say that a subalgebra $L < V$ is a \textbf{Levi subalgebra} if $V =
L \oplus \Rad V$ as vector spaces.  Ling showed that, as in the theory
of Lie algebras, Lie $n$-algebras admit a Levi decomposition.

\begin{theorem}[{\cite[Theorem~4.1]{LingSimple}}]\label{th:Levi}
  Let $V$ be a Lie $n$-algebra.  Then $V$ admits a Levi subalgebra.
\end{theorem}

A further result of Ling's which we shall need is the following.  Let
us say that a Lie $n$-algebra is \textbf{reductive} if its radical
coincides with its centre: $\Rad V = Z(V)$.

\begin{theorem}[{\cite[Theorem~2.10]{LingSimple}}]\label{th:emissingling}
  Let $V$ be a Lie $n$-algebra.  Then $V$ is reductive if and only if
  the Lie algebra $\ad V$ of inner derivations is semisimple.  If in
  addition $\Der V = \ad V$, $V$ is semisimple.
\end{theorem}




In turn this allows us to prove the following useful result.

\begin{proposition}\label{pr:exactness}
  Let $0 \to A \to B \to \Cbar \to 0$ be an exact sequence of Lie
  $n$-algebras.  If $A$ and $\Cbar$ are semisimple, then so is $B$.
\end{proposition}

\begin{proof}
  Since $A$ is semisimple, Theorem~\ref{th:emissingling} says that
  $\ad A$ is semisimple.  $B$ is a representation of $\ad A$, hence
  fully reducible.  Since $A$ is an $\ad A$-submodule of $B$, we have
  $B = A \oplus C$, where $C$ is a complementary $\ad A$-submodule.
  Since $A\lhd B$ is an ideal (being the kernel of a homomorphism),
  $\ad A (C) = 0$, whence $[A \dots A C]=0$.

  The subspace $C$ is actually a subalgebra, since the component
  $[C \dots C]_A$ of $[C \dots C]$ along $A$ is $\ad A$-invariant by
  the $n$-Jacobi identity and the fact that $C$ is $\ad A$-invariant.
  This means that $[C \dots C]_A$ is central in $A$, but $A$ is
  semisimple, whence it must vanish.  Hence, $[C \dots C]\subset C$.
  Since the projection $B \to \Cbar$ maps $C$ isomorphically to
  $\Cbar$, we see that this isomorphism is one of Lie $n$-algebras,
  hence $C<B$ is semisimple and indeed $[C \dots C]=C$.

  Next we show that $[AC\dots C]=0$.  Indeed, for $c_1,\dots,c_{n-1}
  \in C$, the map $a \mapsto [c_1\dots c_{n-1} a]$ is a derivation of
  $A$.  Since $A$ is semisimple, it is an inner derivation.  However
  since $\ad A$ acts trivially on $C$, this derivation is $\ad
  A$-invariant, which means that it is central.  Since $\ad A$ has
  trivial centre, we see that it must be zero.  This shows that $B = A
  \oplus C$ is also a direct sum of $\ad C$-modules, with $A$ being a
  trivial $\ad C$-module.

  Now consider $W_k:=[\underbrace{A \dots A}_{n-k} \underbrace{C \dots
    C}_k]$.  We have seen that $W_0=A$, $W_1 = 0 = W_{n-1}$ and
  $W_n=C$.  We claim that $W_{1<k<n-1}=0$ as well.  Indeed, $W_k$ is
  the image of $\Lambda^{n-k}A \otimes \Lambda^k C \to A$ (since $A$ is
  an ideal) under the bracket.  Since the bracket is $\ad
  V$-equivariant, it is in particular $\ad C$-equivariant, but now $A$
  is a trivial $\ad C$-module and $C$, being semisimple, decomposes
  into nontrivial irreducible $\ad C$-modules.  Therefore the only
  $\ad C$-equivariant map $\Lambda^{n-k}A \otimes \Lambda^k C \to A$,
  for $k\geq 1$, is the zero map.

  In other words, $[ACB\dots B]=0$, whence $B = A \oplus C$ is the direct
  sum of the two commuting ideals $A$ and $C$.  Since $A$ and $C$ are
  themselves direct sum of simple ideals, so is $B$.
\end{proof}

A useful notion that we will need is that of a representation of a Lie
$n$-algebra.  A \textbf{representation} of Lie
$n$-algebra $V$ on a vector space $W$ is a Lie $n$-algebra structure
on the direct sum $V \oplus W$ satisfying the following three
properties:
\begin{enumerate}
\item the natural embedding $V \to V \oplus W$ sending $v \mapsto
  (v,0)$ is a Lie $n$-algebra homomorphism, so that $[V \dots V]
  \subset V$ is the original $n$-bracket on $V$;
\item $[V \dots V W] \subset W$; and
\item $[V \dots V W W]=0$.
\end{enumerate}
We will often say that $W$ is a \textbf{$V$-module}, although this is
slightly misleading in the absence of a notion of a ``universal
enveloping algebra'' for a Lie $n$-algebra.  The second of the above
conditions says that if $W$ is a representation of $V$, we have a map
$\ad V \to \End W$ from inner derivations of $V$ to linear
transformations on $W$.  The $n$-Jacobi identity for $V\oplus W$ says
that this map is a representation of the Lie algebra $\ad V$.
Viceversa, any representation $\ad V \to \End W$ defines a Lie
$n$-algebra structure on $V \oplus W$ extending the Lie $n$-algebra
structure of $V$ and demanding that $[V,\dots,V,W,W]=0$.  Taking $W=V$
gives rise to the \textbf{adjoint representation}, whereas taking $W =
V^*$ gives rise to the \textbf{coadjoint representation}, where if
$\alpha \in V^*$ then
\begin{equation}
  \label{eq:coadjoint}
  [v_1,\dots,v_{n-1},\alpha] = \beta \in V^* \qquad\text{where}\qquad
  \beta(v) = - \alpha\left([v_1,\dots,v_{n-1},v]\right)~.
\end{equation}

\subsection{Basic notions about metric Lie $n$-algebras}
\label{sec:basic-metric}

Let us now introduce an inner product, so that $(V,\Phi,b)$ is a
metric Lie $n$-algebra.

If $W \subset V$ is any subspace, we define
\begin{equation*}
  W^\perp = \left\{v \in V\middle | \left<v,w\right>=0~,\forall w\in
      W\right\}~.
\end{equation*}
Notice that $(W^\perp)^\perp = W$.  We say that $W$ is
\textbf{nondegenerate}, if $W \cap W^\perp = 0$, whence $V = W \oplus
W^\perp$; \textbf{isotropic}, if $W \subset W^\perp$; and
\textbf{coisotropic}, if $W \supset W^\perp$.  Of course, in
positive-definite signature, all subspaces are nondegenerate.

An equivalent criterion for decomposability is the existence of a
proper nondegenerate ideal: for if $I\lhd V$ is nondegenerate, $V = I
\oplus I^\perp$ is an orthogonal direct sum of ideals.  The proofs of
the following results can be read off \emph{mutatis mutandis} from
the similar results for metric Lie $3$-algebras in
\cite[§2.2]{Lor3Lie}.

\begin{lemma}\label{le:coisoquot}
  Let $I\lhd V$ be a coisotropic ideal of a metric Lie $n$-algebra.
  Then $I/I^\perp$ is a metric Lie $n$-algebra.
\end{lemma}

\begin{lemma}\label{le:centreperp}
  Let $V$ be a metric Lie $n$-algebra.  Then the centre is the
  orthogonal subspace to the derived ideal; that is,
  $[V,\dots,V]=Z^\perp$.
\end{lemma}

\begin{proposition}\label{pr:ideals}
  Let $V$ be a metric Lie $n$-algebra and $I \lhd V$ be an ideal.
  Then
  \begin{enumerate}
  \item $I^\perp \lhd V$ is also an ideal;
  \item $I^\perp\lhd Z(I)$; and
  \item if $I$ is minimal then $I^\perp$ is maximal.
  \end{enumerate}
\end{proposition}

\section{Structure of metric Lie $n$-algebras}
\label{sec:structure}

We now investigate the structure of metric Lie $n$-algebras.  If a Lie
$n$-algebra $V$ is not simple or one-dimensional, then it has a proper
ideal and hence a minimal ideal.  Let $I\lhd V$ be a minimal ideal of
a metric Lie $n$-algebra.  Then $I \cap I^\perp$, being an ideal
contained in $I$, is either $0$ or $I$.  In other words, minimal
ideals are either nondegenerate or isotropic.  If nondegenerate, $V =
I \oplus I^\perp$ is decomposable.  Therefore if $V$ is
indecomposable, $I$ is isotropic.  Moreover, by
Proposition~\ref{pr:ideals} (2), $I$ is abelian and furthermore,
because $I$ is isotropic, $[I I V \dots V]=0$.

It follows that if $V$ is euclidean and indecomposable, it is either
one-dimensional or simple, whence of the form \eqref{eq:simple-n-Lie}
with all $\varepsilon_i=1$.  This result, originally due to Nagy
\cite{NagykLie} (see also \cite{GPkLie}), was conjectured in
\cite{FOPPluecker}, albeit in the guise of a conjectural
generalisation of the classical Plücker identity.

Let $V$ be an indecomposable metric Lie $n$-algebra.  Then $V$ is
either simple, one-dimensional (provided the index of the inner
product is $<2$) or possesses an isotropic proper minimal ideal $I$
which obeys $[I I V \dots V]=0$.  The perpendicular ideal $I^\perp$ is
maximal and hence by Lemma~\ref{le:simplequot}, $\Ubar := V/I^\perp$
is simple or one-dimensional, whereas by Lemma~\ref{le:coisoquot}, $\Wbar
:=I^\perp/I$ is a metric Lie $n$-algebra.  The inner product on $V$
induces a nondegenerate pairing $g: \Ubar \otimes I \to \RR$.  Indeed,
let $[u] = u + I^\perp \in \Ubar$ and $v\in I$.  Then we define
$g([u],v) = \left<u,v\right>$, which is clearly independent of the
coset representative for $[u]$.  In particular, $I \cong \Ubar^*$ is
either one- or ($n+1$)-dimensional.  If the signature of the metric of
$\Wbar$ is $(p,q)$, that of $V$ is $(p+r,q+r)$ where $r = \dim I = \dim
\Ubar$.

There are two possibilities for $\Ubar$: either it is one-dimensional or
else it is simple.  We will treat both cases separately.

\subsection{$\Ubar$ is one-dimensional}
\label{sec:U1diml}

If the quotient Lie $n$-algebra $\Ubar=V/I^\perp$ is one-dimensional,
so is the minimal ideal $I$.  Let $u \in V$ be such that $u \not\in
I^\perp$, whence its image in $\Ubar$ generates it.  Because $I \cong
\Ubar^*$ is induced by the inner product, there is $v \in I$ such that
$\left<u,v\right> = 1$.  The subspace spanned by $u$ and $v$ is
therefore nondegenerate, and hence as a vector space we have an
orthogonal decomposition $V = \RR(u,v) \oplus W$, where $W$ is the
perpendicular complement of $\RR(u,v)$.  It is clear that $W \subset
I^\perp$, and that $I^\perp = I \oplus W$ as a vector space.  Indeed,
the projection $I^\perp \to \Wbar$ maps $W$ isomorphically onto
$\Wbar$.

From Proposition~\ref{pr:ideals} (2), it is immediate that
$[u,v,x_1\dots,x_{n-2}]=0=[v,x_1\dots,x_{n-1}]$, for all $x_i \in
W$, whence $v$ is central.  Metricity then implies that the only
nonzero $n$-brackets take the form
\begin{equation}
  \label{eq:n-brackets}
  \begin{aligned}[m]
    [u x_1 \dots x_{n-1}] &= [x_1 \dots x_{n-1}]\\
    [x_1 \dots x_n] &= (-1)^n \left<[x_1 \dots x_{n-1}],x_n\right> v +
    [x_1 \dots  x_n]_W~,
  \end{aligned}
\end{equation}
which defines $[x_1 \dots x_{n-1}]$ and $[x_1 \dots x_n]_W$ and where
$x_i\in W$.  The $n$-Jacobi identity is equivalent to the following
two conditions:
\begin{enumerate}
\item $[x_1 \dots x_{n-1}]$ defines a Lie ($n-1$)-algebra structure
  on $W$, which leaves the inner product invariant due to the
  skewsymmetry of $\left<[x_1 \dots x_{n-1}],x_n\right>$; and
\item $[x_1 \dots x_n]_W$ defines a metric Lie $n$-algebra
  structure on $W$ which is invariant under the ($n-1$)-algebra structure.
\end{enumerate}

As we will see below, this makes $V$ into the \emph{double
  extension} of the metric Lie $n$-algebra $W$ by the one-dimensional
Lie $n$-algebra $\Ubar$.

\subsection{$\Ubar$ is simple}
\label{sec:u-simple}

Consider $I^\perp$ as a Lie $n$-algebra in its own right and let
$R=\Rad I^\perp$ denote its radical.  By Theorem~\ref{th:Levi},
$I^\perp$ admits a Levi subalgebra $L<I^\perp$.  Since $I^\perp\lhd V$
and $R\lhd I^\perp$ is a characteristic ideal, $R \lhd V$.  Indeed,
for all $x_i\in V$, $\ad_{x_1 \dots x_{n-1}}$ is a derivation of
$I^\perp$ (since $I^\perp \lhd V$) and by
Theorem~\ref{th:RadChar}, it preserves $R$.  Let $M = V/R$.
Notice that
\begin{equation*}
  \Ubar = V/I^\perp \cong (V/R)/(I^\perp/R) = M/L~,
\end{equation*}
by the standard homomorphism theorems.  Since $L$ and $\Ubar$ are
semisimple, Proposition~\ref{pr:exactness} says that so is $M$ and
moreover that $M \cong L \oplus \Ubar$.  This means that $R$ is also
the radical of $V$, whence $M$ is a Levi factor of $V$.  This
discussion is summarised by the following commutative diagram with
exact rows and columns:
\begin{equation*}
  \begin{CD}
    @.    0     @.    0     @.    @.  \\
    @.  @VVV         @VVV    @.    @.  \\
    @.    R    @=     R      @.    @.  \\
    @.  @VVV         @VVV    @.    @.  \\
 0 @>>>  I^\perp @>>> V  @>>> \Ubar @>>> 0\\
 @.      @VVV        @VVV      @|   @.\\
   0 @>>>  L    @>>> M  @>>> \Ubar @>>> 0\\
    @.  @VVV         @VVV    @.    @.  \\
    @.    0     @.    0     @.    @.
  \end{CD}
\end{equation*}

The map $M \to \Ubar$ admits a section, so that $M$ has a subalgebra
$\widetilde U$ isomorphic to $\Ubar$ and such that $M = \widetilde U
\oplus L$.  Then the vertical map $V \to M$ also admits a section,
whence there is a subalgebra $U < V$ isomorphic to $\Ubar$ such that
$V = I^\perp \oplus U$ (as vector space).  Furthermore, the inner
product on $V$ pairs $I$ and $U$ nondegenerately, whence $I \oplus U$
is a nondegenerate subspace.  Let $W$ denote its perpendicular
complement, whence $V = W \oplus I \oplus U$.  Clearly $I^\perp = W
\oplus I$, whence the canonical projection $I^\perp \to \Wbar$ maps
$W$ isomorphically onto $\Wbar$.

Let us now write the possible $n$-brackets for $V = W \oplus I \oplus
U$.  First of all, by Proposition~\ref{pr:ideals} (2),
$[V,\dots,V,I^\perp,I]=0$.  Since $U<V$, $[U,\dots,U]\subset U$ and
since $I$ is an ideal, $[U,\dots,U,I]\subset I$.  Similarly, since $W
\subset I^\perp$ and $I^\perp \lhd V$ is an ideal,  $[W,\dots,W]\subset W
\oplus I$.  We write this as
\begin{equation*}
  [w_1 \dots w_n] := [w_1 \dots w_n]_W + \varphi(w_1 \dots w_n)~,
\end{equation*}
where $[w_1 \dots w_n]_W$ defines an $n$-bracket on $W$, which is
isomorphic to the Lie $n$-bracket of $\Wbar = I^\perp/I$, and
$\varphi:\Lambda^nW \to I$ is to be understood as an abelian
extension.  At the other extreme we have the bracket $[U\dots UW] \subset
W$ which makes $W$ into an $\ad U$-module.  Metricity forbids a
nonzero $I$-component  to the above bracket:
\begin{equation*}
  \left<[U\dots UW],U\right> = \left<[U\dots U], W\right> = 0~,
\end{equation*}
since $U$ is a subalgebra.  Finally we have a sequence of brackets
\begin{equation*}
  V_k := [\underbrace{U \dots U}_k \underbrace{W \dots W}_{n-k}]
  \subset W \oplus I~,
\end{equation*}
for $0<k<n-1$.  We notice that
\begin{equation*}
  \left<V_k,U\right> =
  \left<[\underbrace{U \dots U}_k \underbrace{W \dots W}_{n-k}],U\right> =
 \left<[\underbrace{U \dots U}_{k+1} \underbrace{W \dots W}_{n-k-1}],W\right> =
  \left<V_{k+1},W\right>~,
\end{equation*}
whence the component of $V_k$ along $I$ agrees up to a sign with the
component of $V_{k+1}$ along $W$.  In principle all such brackets
occur and the only conditions apart from the metricity come from the
Jacobi identity of $V$.

Similarly to the case when $\Ubar$ is one-dimensional, we will
interpret $V$ as the \emph{double extension} of the metric Lie
$n$-algebra $\Wbar$ by the simple Lie $n$-algebra $U$.

\subsection{Double extensions and the structure theorem}
\label{sec:doublext}

More generally we have the following definition.
\begin{definition}
  Let $W$ be a metric Lie $n$-algebra and let $U$ be a Lie
  $n$-algebra.  Then by the \textbf{double extension of $W$ by $U$} we
  mean the metric Lie $n$-algebra on the vector space $W \oplus U
  \oplus U^*$ with the following nonzero $n$-brackets subject to the
  Jacobi identity for $V$:
  \begin{itemize}
  \item $[U\dots U] = [U\dots U]_U$, making $U$ into a subalgebra;
  \item $[U\dots U U^*] \subset U^*$, making $U^*$ into the coadjoint
    representation of $U$;
  \item $[U\dots U W] \subset W$, making $W$ into an $\ad U$-module;
  \item $[w_1\dots w_n] = [w_1\dots w_n]_W + \varphi(w_1,\dots,w_n)$
    for $w_i\in W$, where $[\dots]_W$ is the bracket of the Lie
    $n$-algebra $W$ and $\varphi: \Lambda^n W \to U^*$ is an $\ad
    U$-equivariant map; and
  \item $\ad U$-equivariant brackets
    $[\underbrace{U \dots U}_k \underbrace{W \dots W}_{n-k}] \subset W
    \oplus U^*$ for $0<k<n-1$, where metricity identifies (perhaps up
    to a sign) the $W$ component of
    $[\underbrace{U \dots U}_k \underbrace{W \dots W}_{n-k}]$ with the
    $U^*$ component of
    $[\underbrace{U \dots U}_{k-1} \underbrace{W \dots W}_{n-k+1}]$.
  \end{itemize}
  The resulting Lie $n$-algebra is metric, with inner product which
  extends the one on $W$ by the dual pairing between $U$ and $U^*$.
  One can also add any invariant symmetric bilinear form on $U$, even
  if degenerate.
\end{definition}

For $n=2$ this construction is due to Medina and Revoy
\cite{MedinaRevoy}, whereas for $n=3$ it is due to the authors of
\cite{2p3Lie}.

In summary we have proved the following

\begin{theorem}\label{th:metricNLie}
  Every indecomposable metric Lie $n$-algebra is either
  one-dimensional, simple or else it is the double extension of a
  metric Lie $n$-algebra of smaller dimension by a one-dimensional or
  a simple Lie $n$-algebra.
\end{theorem}

An easy induction argument and the fact that metric Lie $n$-algebras
are orthogonal direct sums of their indecomposable components yields
the following

\begin{corollary}\label{co:metric3lie}
  The class of metric Lie $n$-algebras is generated by the simple and
  one-dimensional Lie $n$-algebras under the operations of orthogonal
  direct sum and double extension.
\end{corollary}

It is clear that the subclass of euclidean metric Lie $n$-algebras is
generated by the simple and one-dimensional euclidean Lie 3-algebras
under orthogonal direct sum, since double extension always incurs in
indefinite signature.  Therefore an indecomposable euclidean metric Lie
$n$-algebra is either one-dimensional or simple
\cite{FOPPluecker,NagykLie,GPkLie}.  The lorentzian indecomposables
admit at most one double extension by a one-dimensional Lie
$n$-algebra and are easy to classify \cite{Lor3Lie,JMFLorNLie}.

\bibliographystyle{utphys}
\bibliography{AdS,AdS3,ESYM,Sugra,Algebra,Geometry}

\end{document}